\documentclass[12pt,leqno]{article}

\usepackage{amsthm}

\usepackage{amsmath}
\usepackage{amscd}
\usepackage{amsmath}
\usepackage{amssymb}
\usepackage{latexsym}
\makeatletter

\@addtoreset{equation}{section}
\makeatother

\newtheorem{thm}{Theorem}[section]
\newtheorem{pr}[thm]{Proposition}

\newtheorem{lm}[thm]{Lemma}
\newtheorem{cor}[thm]{Corollary}

\newcommand{\Jac}{{\rm Pic}_C}
\newcommand{\Pic}{{\rm Pic}_{C,\mathfrak{m}}}
\newcommand{\Sch}{{\it Sch}}

\input xy \xyoption{all} \CompileMatrices

\begin{document}

\title{Blow-ups and the class field theory for curves}
\author{Daichi Takeuchi}
\date{}
\maketitle

\begin{abstract}
We propose another proof of the geometric class field theory for curves by considering blow-ups 
of  symmetric products of curves.
\end{abstract}

\section{Introduction}

Let $k$ be a perfect field, and $C$ be a projective smooth geometrically connected curve over $k$. The geometric class field theory gives a geometric description of abelian coverings of 
$C$ by using generalized jacobian varieties. Let us recall its precise statement. Fix a modulus $\mathfrak{m}$, i.e. an effective Cartier divisor of $C$ and let $U$ be its complement in $C$. Let $\Pic^0$ be the corresponding 
generalized jacobian variety. Let $G^0\to\Pic^0$ be an \'etale isogeny of smooth commutative algebraic groups and $G^1\to\Pic^1$ be a compatible morphism of torsors. We call such a pair $(G^0\to\Pic^0,G^1\to\Pic^1)$ a {\it covering of }$(\Pic^0,\Pic^1)$. We call such a covering {\it connected abelian} if $G^0$ is 
connected and $G^0\to\Pic^0$ is an abelian isogeny. There is a natural map from $U$ to $\Pic^1$ sending a point of $U$ to 
its associated invertible sheaf with a trivialization. The geometric class field theory states:
\begin{thm}\label{main}
Let $C$ be a projective smooth geometrically connected curve over a perfect field $k$. Fix a modulus $\mathfrak{m}$ of $C$ and denote its complement by $U$. Let $\Pic^0$ be the generalized jacobian variety with modulus 
$\mathfrak{m}$. Then a connected abelian covering $(G^0\to\Pic^0,G^1\to\Pic^1)$ pulls back by the natural map $U\to\Pic^1$ to a geometrically connected abelian covering of $U$ whose ramification is 
bounded by $\mathfrak{m}$. Conversely, every such covering is obtained in this way.
\end{thm}

Originally this theorem was proved by M. Rosenlicht \cite{Ros}. S. Lang \cite{Lang} generalized his results to an 
arbitrary algebraic variety. Their works are explained in detail in Serre's book \cite{Serre}.

On the other hand, in 1980s, P. Deligne found another proof for the tamely ramified case by using symmetric powers of curves. 
The aim of this paper is to complete his proof by considering blow-ups 
of symmetric powers of curves. 

Recently Q. Guignard
 did a similar work to this paper, although we do not know 
his results in detail.

Actually we prove a variant of Theorem \ref{main} now stated. 
\begin{thm}\label{main2}
There is an isomorphism of groups between the subgroup of ${\rm H}^1(U,\mathbb{Q}/\mathbb{Z})$ 
consisting of a character $\chi$ such that $
{\rm Sw}_{P}(\chi)\leq n_P -1$ for all points $P\in \mathfrak{m}$, where $n_P$ is the multiplicity of $\mathfrak{m}$ at $P$, and the subgroup of ${\rm H}^1(\Pic,\mathbb{Q}/\mathbb{Z})$ consisting of $\rho$ which is multiplicative, i.e. the self-external product $\rho\boxtimes 1+ 1
\boxtimes\rho$ on $\Pic\times_k\Pic$ equals to  $m^\ast\rho$,  the pull back of $\rho$ by the multiplication map $m:\Pic\times_k\Pic\to\Pic$.
\end{thm}

The relation between Theorem \ref{main} and Theorem \ref{main2} will be explained in Section \ref{pro}. 

When $k$ is algebraically closed, Theorem \ref{main2} can be stated as follows. 
Let $\rho$ be a multiplicative character of $\Pic$. 
Fix a closed point $P\in\Pic^1$. The multiplicativity of $\rho$ implies that, for an integer $d$, the pull back 
of $\rho^d$ by the  multiplication by $dP$ $\Pic^0\to\Pic^d$ coincides with $\rho^0$. 
In this way, Theorem \ref{main2} can be restated as below:
\begin{thm}\label{main3}
Assume that $k$ is algebraically closed. Then there is an isomorphism of groups between the subgroup of ${\rm H}^1(U,\mathbb{Q}/\mathbb{Z})$ 
consisting of a character $\chi$ such that $
{\rm Sw}_{P}(\chi)\leq n_P -1$ for all points $P\in \mathfrak{m}$ and the subgroup of ${\rm H}^1(\Pic^0,\mathbb{Q}/\mathbb{Z})$ consisting of a multiplicative character $\rho^0$, i.e. 
the self-external product $\rho^0\boxtimes 1+1\boxtimes\rho^0$ on $\Pic^0\times_k\Pic^0$ equals to  $m^\ast\rho^0$, the pull back of $\rho^0$ by the multiplication map $m:\Pic^0\times_k\Pic^0\to\Pic^0$. 
\end{thm}

Throughout this paper, we use the following convention.
 We identify an effective Cartier divisor with the associated closed subscheme. For an object defined on a scheme $S$ (e.g. an $S$-scheme, a locally free sheaf, a vector bundle, and so on) 
 and a $S$-scheme $T$, we denote its pull back to $T$ by the same letter, unless 
 there may be ambiguity.
 We denote 
 the category of $S$-schemes by $\Sch/S$. For a category $\mathcal{C}$, we call a functor 
 $\mathcal{C}^{{\it op}}\to ({\it Set})$, from the opposite category of $\mathcal{C}$ to 
 the category of sets $({\it Set})$, a presheaf on $\mathcal{C}$.
 \section{Preliminaries}\label{wit}

  In this section, we recall the notion of Witt vectors and refined Swan conductors, and prove some formulas. Fix a prime number $p$.

 Let $(R,\pi)$ be a DVR of equal characteristic $p$ and $K$ be its field of fractions. Let $v_R$ be its 
 normalized valuation. 
 Let $m$ be an integer $\geq 0$. We extend the valuation $v_R$ to ${\rm W}_{m+1}(K)$ by setting 
 \begin{equation*}
 v_R((a_0,\dots,a_m)):=\underset{i}\min\{p^{m-i}v_R(a_i)\}.
 \end{equation*}
 We define an increasing exaustive filtration on ${\rm W}_{m+1}(K)$ 
 by setting, for $n\in\mathbb{Z}$, ${\rm fil}_n{\rm W}_{m+1}(K)$ to be the subgroup of ${\rm W}_{m+1}(K)$ consisting of an element $(a_0,\dots,a_m)$ such that 
 \begin{equation*}
 v_R((a_0,\dots,a_m))\geq -n.
 \end{equation*}
 Let $\mathcal{O}_K$ be the structure sheaf of rings on the \'etale topos of ${\rm Spec}(K)$. Let 
 $F$ be the absolute Frobenius map $\mathcal{O}_K\to \mathcal{O}_K$, i.e. sending $x\mapsto x^p$, and let the ring homomorphism ${\rm W}_{m+1}(\mathcal{O}_K)\to{\rm W}_{m+1}(\mathcal{O}_K)$ induced from $F$ denoted by the same letter $F$.
  The short exact sequence 
  \begin{equation*}
  0\to\mathbb{Z}/p^{m+1}\mathbb{Z}\to{\rm W}_{m+1}(\mathcal{O}_K)\overset{F -1}\to{\rm W}_{m+1}(\mathcal{O}_K)\to 0
  \end{equation*}
  of \'etale sheaves on ${\rm Spec}(K)$ defines an isomorphism ${\rm W}_{m+1}(K)/{\rm Im}(F-1)
  \overset\cong\to {\rm H}^1(K,\mathbb{Z}/p^{m+1}\mathbb{Z})$. Define an increasing exaustive 
  filtration ${\rm fil}_n{\rm H}^1(K,\mathbb{Z}/p^{m+1}\mathbb{Z})$ of ${\rm H}^1(K,\mathbb{Z}/p^{m+1}\mathbb{Z})$ by the image of ${\rm fil}_n{\rm W}_{m+1}(K)$ through the boundary map $\delta_{m+1,K}$.
  
  For any $\chi\in{\rm H}^1(K,\mathbb{Z}/p^{m+1}\mathbb{Z})$, the Swan conductor of $\chi$, 
  ${\rm Sw}_R(\chi)$, is the smallest integer $n\geq 0$ such that $\chi\in{\rm fil}_n{\rm H}^1(K,\mathbb{Z}/p^{m+1}\mathbb{Z})$(cf. \cite{Kato}). When $R$ is henselian and the residue field is perfect, this is the same as the classical Swan conductor.
  \begin{lm}\label{fil}
  Let $R$ and $K$ be as above. Take $\chi\in{\rm H}^1(K,\mathbb{Z}/p^{m+1}\mathbb{Z})$.
  
  1. Let $\widehat{R}$ be the completion of $R$ and $\widehat{K}$ be its field of fractions. Denote the restriction of 
  $\chi$ to $\widehat{K}$ by $\widehat{\chi}$. Then, the equality ${\rm Sw}_R(\chi) ={\rm Sw}_{\widehat{R}}(\widehat{\chi})$ holds.
  
  2. The subgroup ${\rm fil}_0{\rm H}^1(K,\mathbb{Z}/p^{m+1}\mathbb{Z})$ of ${\rm H}^1(K,\mathbb{Z}/p^{m+1}\mathbb{Z})$ coincides with the image of 
  the map ${\rm H}^1({\rm Spec}(R),\mathbb{Z}/p^{m+1}\mathbb{Z})\to {\rm H}^1(K,\mathbb{Z}/p^{m+1}\mathbb{Z})$, i.e. the group of unramified characters.
  \end{lm}
  \proof{
  
  1.The commutative diagram 
  \begin{equation*}
  \xymatrix{
  {\rm W}_{m+1}(K)\ar[r]^{\delta_{m+1,K}\ \ \ }\ar[d]&{\rm H}^1(K,\mathbb{Z}/p^{m+1}\mathbb{Z})\ar[d]\\
    {\rm W}_{m+1}(\widehat{K})\ar[r]^{\delta_{m+1,\widehat{K}}\ \ \ }&{\rm H}^1(\widehat{K},\mathbb{Z}/p^{m+1}\mathbb{Z})
    }
    \end{equation*}
    implies ${\rm Sw}_R(\chi)\geq{\rm Sw}_{\widehat{R}}(\widehat{\chi})$. Let $n={\rm Sw}
    _{\widehat{R}}(\widehat{\chi})$. Then there exists a Witt 
    vector $\widehat{\alpha}\in {\rm fil}_n{\rm W}_{m+1}(\widehat{K})$ mapping to $\widehat{\chi}$. Take $\alpha\in {\rm fil}_n{\rm W}_{m+1}(K)$ so that every component of $\alpha$ is close enough to that of $\widehat{\alpha}$.
If $\widehat{\alpha}-\alpha$ (here $\alpha$ is regarded as an element of ${\rm W}_{m+1}(\widehat{K})$) is in ${\rm W}_{m+1}(\widehat{R})$, $\delta_{m+1,\widehat{K}}(\widehat{\alpha}-\alpha)$ is an unramified character by 2. Therefore, $\chi-\delta_{m+1,K}(\alpha)$ is unramified. Again by 2., there exists $\beta\in{\rm W}_{m+1}(R)$ such that $\chi-\delta_{m+1,K}(\alpha)=\delta_{m+1,K}(\beta)$, hence the assertion.

2. This follows from the commutative diagram
\begin{equation}\label{c}
\xymatrix{
{\rm W}_{m+1}(R)\ar[r]\ar[d]&{\rm H}^1({\rm Spec}(R),\mathbb{Z}/p^{m+1}\mathbb{Z})\ar[d]\\
{\rm W}_{m+1}(K)\ar[r]&{\rm H}^1(K,\mathbb{Z}/p^{m+1}\mathbb{Z})
}
\end{equation}
and the fact that the two horizontal arrows in (\ref{c}) are surjective.
\qed}

Denote $\widehat{\Omega}^{1}_R$ to be 
the $\pi$-adic completion of the absolute differential module $\Omega^1_R$. Let $\widehat{\Omega}^1_K:=
\widehat{\Omega}^{1}_R\otimes_R K$. 
The canonical map $\widehat{\Omega}^1_R\to\widehat{\Omega}^1_K$ is injective 
and we usually regard $\widehat{\Omega}^1_R$ as an $R$-submodule of 
$\widehat{\Omega}^1_K$ via this map. The $R$-module $\widehat{\Omega}^{1}_R({\rm log})$ is the 
$R$-submodule of $\widehat{\Omega}^{1}_K$ generated by $\widehat{\Omega}^{1}_R$ and 
$d{\rm log}\pi:=\frac{d\pi}{\pi}$. From the definition, the following holds:
\begin{lm}\label{diflog}
Assume that $R$ is obtained from a smooth scheme over a perfect field by localizing at a 
point of codimension one. Let $b_1,\dots,b_n$ be a lift of a $p$-basis of the residue field of $R$ to $R$. Then, 
$\widehat{\Omega}^1_R({\rm log})$ is a $\widehat{R}$-free module with a basis $db_1,\dots,
db_n,d{\rm log}\pi$.
\end{lm}

\qed

For $\omega\in\widehat{\Omega}^{1}_K$, define $v^{\rm log}_R(\omega)$ as the 
largest integer $n$ such that $\omega\in \pi^n\widehat{\Omega}^{1}_R({\rm log})$ (we formally 
put $v_R^{{\rm log}}(0):=\infty$). There is a homomorphism $F^m d:{\rm W}_{m+1}(K)\to
\widehat{\Omega}^{1}_K$ given by 
\begin{equation*}
F^m d((a_0,\dots,a_m)):=\sum_i a_i^{p^{m-i}-1}da_i.
\end{equation*}
Define an increasing exaustive filtration on $\widehat{\Omega}^1_K$ by setting 
\begin{equation*}
{\rm fil}_n\widehat{\Omega}^1_K:=\{\omega\in\widehat{\Omega}^1_K\lvert 
v^{\rm log}_R(\omega)\geq-n\}
\end{equation*}
for $n\in\mathbb{Z}$.
The homomorphism $F^m d:{\rm W}_{m+1}(K)\to \widehat{\Omega}^1_K$ respects their filtrations. 
In other words, $v_R(\alpha)\leq v^{\rm log}_R(F^m d\alpha)$ hold for all $\alpha\in{\rm W}_{m+1}(K)$. 

The next aim of this section is to prove the formulas in Corollary \ref{blprod} and \ref{Witt2}. To do 
this, we need the following refinement of the refined Swan conductor in \cite{Kato}, proved in \cite{Leal}.
\begin{pr}\label{Lea}

  Let $n$ be an integer $\geq 0$.

1. There is a unique homomorphism 
\begin{equation*}
{\rm rsw}:{\rm fil}_n{\rm H}^1(K,\mathbb{Z}/p^{m+1}\mathbb{Z})\to
{\rm fil}_n\widehat{\Omega}^1_K /{\rm fil}_{\lfloor\frac{n}{p}\rfloor}\widehat{\Omega}^1_K
\end{equation*}
such that the composition 
\begin{equation*}
{\rm fil}_n{\rm W}_{m+1}(K)\to{\rm fil}_n{\rm H}^1(K,\mathbb{Z}/p^{m+1}\mathbb{Z})\to
{\rm fil}_n\widehat{\Omega}^1_K /{\rm fil}_{\lfloor\frac{n}{p}\rfloor}\widehat{\Omega}^1_K
\end{equation*}
coincides with $F^m d$.

2. For $\lfloor\frac{n}{p}\rfloor\leq i\leq n$, the induced map 
\begin{equation*}
{\rm fil}_n{\rm H}^1(K,\mathbb{Z}/p^{m+1}\mathbb{Z})/{\rm fil}_i{\rm H}^1(K,\mathbb{Z}/p^{m+1}\mathbb{Z})\to
{\rm fil}_n\widehat{\Omega}^1_K/{\rm fil}_i\widehat{\Omega}^1_K
\end{equation*}
is injective.
\end{pr}
Let $k$ be a perfect field of characteristic $p$. 
Let $X,Y$ be smooth schemes over $k$. Let $D$ (resp. $E$) be a 
smooth irreducible closed subvariety of $X$ (resp. $Y$). Let $\widetilde{X}$ (resp. $\widetilde{Y}$, resp. 
$\widetilde{X\times Y}$) be the blow-up of $X$ (resp. $Y$, resp. $X\times Y$) along $D$ (resp. 
$E$, resp. $D\times E$). Let $R_1$ (resp. $R_2$, resp. $R_3$) be the DVR at the generic point 
of the exceptional divisor of $\widetilde{X}$ (resp. $\widetilde{Y}$, resp. $\widetilde{X\times Y}$). 
Let $K_i$ be the field of fractions of $R_i$ for $i=1,2,3$.

\begin{lm}\label{dprod}
The projections $X\times Y\to X$ and $X\times Y\to Y$ induce an isomorphism 
\begin{equation*}
\widehat{\Omega}^1_{K_3}\cong (\widehat{K_3}\otimes_{\widehat{K_1}}
\widehat{\Omega}^1_{K_1})\oplus (\widehat{K_3}\otimes_{\widehat{K_2}}\widehat{\Omega}^1_{K_2}).
\end{equation*}
This isomorphism 
respects the filtrations, i.e. via this isomorphism, ${\rm fil}_n\widehat{\Omega}^1
_{K_3}$ coincides with $(\widehat{R}_3\otimes_{\widehat{R}_1}{\rm fil}_n\widehat{\Omega}^1_{K_1})
\oplus(\widehat{R}_3\otimes_{\widehat{R}_2}{\rm fil}_n\widehat{\Omega}^1_{K_2})$.
\end{lm}
\proof{
The first assertion follows from the isomorphism 
\begin{equation*}
\Omega^1_{X\times Y}\cong{\it pr}^\ast_X\Omega^1_X\oplus{\it pr}^\ast_Y\Omega^1_Y,
\end{equation*}
where ${\it pr}_X$ and ${\it pr}_Y$ are the projections to $X$ and $Y$. 
Note that ${\it pr}_X$ and ${\it pr}_Y$ induce extensions $R_3/R_1$ and $R_3/R_2$ of DVRs , 
which preserve uniformizers. The assertion follows.

\qed}
\begin{cor}\label{blprod}
 Let $\chi_i\in{\rm H}^1(K_i,\mathbb{Q}/
\mathbb{Z})$ for $i=1,2$. Then, the following holds:
\begin{equation*}
{\rm Sw}_{R_3}(\chi_1\boxtimes 1+1\boxtimes\chi_2)=
\max\{{\rm Sw}_{R_1}(\chi_1),{\rm Sw}_{R_2}(\chi_2)\}.
\end{equation*}
\end{cor}
\proof{
 Since the prime-to-$p$ parts of 
$\chi_1,\chi_2$, and $\chi_1\boxtimes 1+1\boxtimes\chi_2$ are tame, 
we reduce to the case when $\chi_i\in{\rm H}^1(K_i,\mathbb{Z}/p^{m+1}\mathbb{Z})$. 
Since the extensions $R_3/R_1,R_3/R_2$ of DVRs preserve uniformizers, the 
inequality ${\rm Sw}_{R_3}(\chi_1\boxtimes 1+1\boxtimes\chi_2)\le\max\{{\rm Sw}_{R_1}(\chi_1),{\rm Sw}_{R_2}(\chi_2)\}$ is obvious.  Proposition 
\ref{Lea} and Lemma \ref{dprod} imply the assertion.

\qed}

Let $S$ be a scheme. For a quasi-projective $S$-scheme $X$ and a natural number $d\geq 1$, the $d$-th symmetry group $\mathfrak{S}_{d}$ 
acts on $X^d:=X\times_S X\times_S\dots\times_S X$ ($d$ times) via permutation of coordinates.
 Define a scheme $X^{(d)}:=X^{d}/\mathfrak{S}_d$. $X^{(d)}$ is called the $d$-th symmetric product of $X$. It is known that, if $X$ is smooth of relative dimension $1$ over $S$, $X^{(d)}$ is smooth and parametrizes effective Cartier divisors of $\deg=d$ on $X$. In particular, the formation of $X^{(d)}$ commutes with base change $S'\to S$.

Let $C$ be a projective smooth geometrically 
 connected curve over $k$.
Let $U$ be a non-empty open subscheme of $C$. Take and fix a $k$-rational point $P$ of $C$ outside $U$.  

Let $d$ be an integer $\geq 1$. 
 We construct a map
  ${\rm H}^{1}(U, \mathbb{Q}/\mathbb{Z})\to{\rm H}^{1}(U^{(d)}, \mathbb{Q}/\mathbb{Z})$ as follows. First fix a finite abelian group $G$.
 Let $V\to U$ be a $G$-torsor.
 Then $V^{d}$ is a $G^{d}$-torsor of $U^{d}$.
 Let $H$ be the subgroup of  $G^{d}$ consisting of elements
  $(a_{1}, \dots ,a_{d})$ satisfying $\sum_{1\leq i\leq d}a_{i}=0$. Then $V^{d}/H$ is a  $G$-torsor of $U^{d}$. This torsor has a natural action by the $d$-th symmetry group $\mathfrak{S}_{d}$ which is equivariant with respect to its action to $U^{d}$. 
\begin{lm}
  The morphism 
  \begin{equation}
  (V^d/H)/\mathfrak{S}_d\to U^{(d)}
  \end{equation}
  induced from the map $V^d/H\to U^d$, taking the quotients by $\mathfrak{S}_d$, is 
  a $G$-torsor. 
\end{lm}
\proof{
It is sufficient to show that, for every geometric point $\bar{x}$ of $U^d$, the inertia group 
$(\mathfrak{S}_d)_{\bar{x}}$ at $\bar{x}$ acts trivially on the fiber $(V^d/H)_{\bar{x}}$ over 
$\bar{x}$, see \cite[Remarque 5.8.]{SGA1}. 

We may assume that $k$ is algebraically closed and that geometric points considered are 
$k$-valued points. Let $\bar{x}$ be a geometric point of $U^d$. For simplicity, we assume that 
$\bar{x}=(x_1,\dots,x_1,x_2,\dots,x_2,\dots,x_r,\dots,x_r)$, where $x_1,\dots,x_r$ are distinct points and 
$x_i$ appears $d_i$ times for each $i$. Then the inertia group $(\mathfrak{S}_d)_{\bar{x}}$ at $\bar{x}$ is isomorphic to $\prod_{1\le i\le r}\mathfrak{S}_{d_i}$. 

For each $i$, take a 
$k$-valued point $e_i$ of $V\times_U x_i$. From the definition of $H$, 
the fiber of $V^d/H$ over $\bar{x}$ can be identified with the set 
\begin{equation}
\left\{(e_1,e_1,\dots,e_r,ge_r)\lvert g\in G\right\},
\end{equation}
on which $(\mathfrak{S}_d)_{\bar{x}}$ acts trivially. 
\qed}
  
  In this way, we construct a $G$-torsor $(V^{d}/ H)/\mathfrak{S}_{d}$ on $U^{(d)}$. 
Since this construction is compatible with a morphism of abelian groups $G\to G'$, we obtain a 
group homomorphism  ${\rm H}^{1}(U, \mathbb{Q}/\mathbb{Z})\to{\rm H}^{1}(U^{(d)}, \mathbb{Q}/\mathbb{Z})$. Let this map denoted by $\kappa$. 
We also 
denote $\chi^{(d)}:=\kappa(\chi)$. 

We consider a similar construction on the groups of Witt vectors. Let $K$ be the field of fractions of $U$, $K'$ be that of $U^{(d)}$, and $K''$ be that of $U^d$. Denote the morphism $K\to K''$, 
induced by the $i$th projection $U^d\to U$, by ${\it pr}_{i}^{\ast}$. Consider the map 
$\lambda:{\rm W}_{m+1}(K)\to {\rm W}_{m+1}(K'')$ sending a Witt vector $\alpha$ to ${\it pr}_1^\ast\alpha+\dots+{\it pr}_{d}^\ast\alpha$. Since the extension $K''/K'$, induced by 
the natural projection $U^d\to U^{(d)}$, is finite Galois with the 
Galois group $\mathfrak{S}_d$, the $\mathfrak{S}_d$-fixed part of ${\rm W}_{m+1}(K'')$ coincides with ${\rm W}_{m+1}(K')$ (here ${\rm W}_{m+1}(K')$ is considered as a subgroup of ${\rm W}_{m+1}(K'')$ via the natural projection $U^{(d)}\to U^d$). Thus 
the map $\lambda$ factors through ${\rm W}_{m+1}(K')$ and commutes with the boundary maps. We also denote the induced map ${\rm W}_{m+1}(K)\to{\rm W}_{m+1}(K')$ by $\lambda$. Also, the 
canonical morphism $\widehat{\Omega}^1_{K'}
\otimes_{K'}K''
\to\widehat{\Omega}^1_{K''}$ is an isomorphism and the $\mathfrak{S}_d$-fixed part of 
$\widehat{\Omega}^1_{K''}$ coincides with (the image of) $\widehat{\Omega}^1_{K'}$.
We define a map $\mu:\widehat{\Omega}^1_K\to\widehat{\Omega}^1_{K'}$ similarly to $\lambda$. 
The maps $\lambda$ and $\mu$ commute with $F^md$.

Let $R$ be the DVR of $C$ at $P$, and $R'$ be the DVR of $K'$ at the generic point of the 
exceptional divisor of the blow-up of $C^{(d)}$ along the point corresponding to the divisor $dP$. We define filtrations on ${\rm W}_{m+1}(K)$ (resp. ${\rm W}_{m+1}(K')$) and $\widehat{\Omega}^1_K$ (resp. $\widehat{\Omega}^1_{K'}$) by $R$ (resp. $R'$) (cf. Section \ref{wit}).

\begin{thm}
\label{Witt}
Let $n$ be an integer. 

1. The homomorphism 
\begin{equation*}
\lambda:{\rm W}_{m+1}(K)\to{\rm W}_{m+1}(K')
\end{equation*}
 sends ${\rm fil}_n{\rm W}_{m+1}(K)$ into ${\rm fil}_{\left\lfloor\frac{n}{d}\right\rfloor}{\rm W}_{m+1}(K')$.

2. The homomorphism 
\begin{equation*}
\mu:\widehat{\Omega}^1_K\to\widehat{\Omega}^1_{K'}
\end{equation*}
 sends  ${\rm fil}_n\widehat{\Omega}^1_K$ into ${\rm fil}_{\left\lfloor\frac{n}{d}\right\rfloor}\widehat{\Omega}^1_{K'}$. Let 
$j$ be an integer. The induced map 
\begin{equation*}
{\rm fil}_{(j+1)d-1}\widehat{\Omega}^1_K/
{\rm fil}_{jd}\widehat{\Omega}^1_K\to{\rm Gr}_{j}\widehat{\Omega}^1_{K'}
\end{equation*} 
is injective, here ${\rm Gr}_j:={\rm fil}_j/{\rm fil}_{j-1}$.
\end{thm}

\begin{cor}\label{Witt2}
Let $\chi$ be a character in ${\rm H}^{1}(U, \mathbb{Q}/\mathbb{Z})$.
The following identity holds:
\begin{equation*}
{\rm Sw}_{R'}(\chi^{(d)})=\left\lfloor\frac{{\rm Sw}_R(\chi)}{d}\right\rfloor.
\end{equation*}

\end{cor}

(Proof of Corollary \ref{Witt2})
Taking the prime-to-$p$ part of $\chi$, we reduce to the case when $\chi\in
{\rm H}^1(U,\mathbb{Z}/p^{m+1}\mathbb{Z})$. 
Take $\alpha\in{\rm W}_{m+1}(K)$ such that $\alpha$ maps to $\chi$ via the boundary map 
${\rm W}_{m+1}(K)\to{\rm H}^1(U,\mathbb{Z}/p^{m+1}\mathbb{Z})$ and 
$v_R(\alpha)=-{\rm Sw}_R(\chi)$. In this case, the equality $v^{{\rm log}}_{R}(F^md\alpha)=
-{\rm Sw}_R(\chi)$ holds by Proposition \ref{Lea}.2. Let $r:=\lfloor\frac{{\rm Sw}_R(\chi)}{d}\rfloor$. 
When ${\rm Sw}_R(\chi)=0$, $\chi$ is unramified. Thus $\chi^{(d)}$ is unramified too by the construction of $\chi^{(d)}$, which imply the assertion. Assume ${\rm Sw}_R(\chi)>0$.

 Consider the following commutative diagram 
\begin{equation}\label{sw}
\xymatrix{
{\rm fil}_{r}{\rm W}_{m+1}(K')\ar[r]\ar[rd]_{F^m d}&{\rm fil}_r{\rm H}^1(K',\mathbb{Z}/p^{m+1}\mathbb{Z})/
{\rm fil}_{r-1}{\rm H}^1(K',\mathbb{Z}/p^{m+1}\mathbb{Z})\ar[d]^{{\rm rsw}}\\
&
{\rm fil}_r\widehat{\Omega}^1_{K'}/{\rm fil}_{r-1}\widehat{\Omega}^1_{K'}.
}
\end{equation}
From Theorem \ref{Witt}.1, the inequality $v_{R'}(\lambda\alpha)\ge -r$ holds. Thus $\chi^{(d)}$ is contained in $\allowbreak{\rm fil}_r{\rm H}^1(K',\mathbb{Z}/p^{m+1}\mathbb{Z})$. Since 
the image of $\chi^{(d)}$ by ${\rm rsw}$ in the diagram (\ref{sw}) coincides with the class containing $F^md\lambda\alpha=\mu F^md\alpha$, the assertion follows from Theorem \ref{Witt}.2. 
and the equality $v^{{\rm log}}_R(F^md\alpha)=-{\rm Sw}_R(\chi)$.

\qed

First we prove some lemmas to prove Theorem \ref{Witt}. Let $R''$ be the 
normalization of $R'$ in $K''$. $R''$ is a DVR. The natural projection $C^d\to C^{(d)}$ 
and the $i$th projection $C^d\to C$ define extensions of DVRs
\begin{equation*}
R'\hookrightarrow R''\overset{{\it pr}^\ast_i}\hookleftarrow R.
\end{equation*}
Fix a uniformizer $t$ of $R$. Let $S_1,\dots,S_d$ be the elementary symmetric polynomials of 
${\it pr}^\ast_1 t,\dots{\it pr}^\ast_d t$ in $R''$, i.e. $S_1,\dots,S_d$ satisfy the following identity
\begin{equation*}
(T-{\it pr}^\ast_1 t)\cdots(T-{\it pr}^\ast_d t)=T^d -S_1 T^{d-1}+\dots+(-1)^d S_d.
\end{equation*}

\begin{lm}\label{unif}
The elements $S_1,\dots,S_d$ are uniformizers of $R'$. 
The valuations of ${\it pr}^\ast_1 t,\dots,{\it pr}^\ast_d 
t$ with respect to $R''$ are the same. 
\end{lm}
\proof{
The first assertion is obvious from the definition of the blow-up. 
The permutation
\begin{equation*}
\sigma:C^d\to C^d
 \end{equation*}
 of the first and the $i$th coordinates induces an isomorphism 
 of extensions of DVRs
 \begin{equation*}
 \xymatrix{
 R''\ar[rr]_{\cong}^{\sigma^\ast}&&R''\\
 &R.\ar[lu]^{{\it pr}^\ast_i}\ar[ru]_{{\it pr}^\ast_1}
 }
 \end{equation*}
 The second assertion follows.
\qed}
 
Note that, by Lemma \ref{diflog} and the fact that the residue field of $R'$ is $k(\frac{S_1}{S_d},
\dots,\frac{S_{d-1}}{S_d})$, $\widehat{\Omega}^1_{R'}({\rm log})$ is an $R'$-free module with a 
basis $\frac{dS_1}{S_d},\dots,\frac{dS_d}{S_d}$. 

\begin{lm}\label{anbasis}
For each integer $i$, define 
\begin{equation*}
\omega_i:=\frac{d({\it pr}^\ast_1 t)}{{{\it pr}^\ast_1 t}^i}+\dots 
+\frac{d({\it pr}^\ast_d t)}{{{\it pr}^\ast_d t}^i}.
\end{equation*}
Let $j$ be an integer. Then, the differentials $\omega_{jd+1},\dots,\omega_{(j+1)d}$ form 
an $R'$-basis of the $R'$-free module $\frac{1}{S_d^j}\widehat{\Omega}^1_{R'}({\rm log})$.

 \end{lm}
 \proof{
 Since the differentials $\omega_j$ are $\mathfrak{S}_d$-invariant, they are indeed contained in 
 $\widehat{\Omega}^1_{K'}$. 
 
Suppose $j\ge0$. Define a polynomial $F(T):=(T-{\it pr}^\ast_1 t)\cdots(T-{\it pr}^\ast_d t)$. The following equalities 
 hold: 
 \begin{align*}
 -dS_1 T^{d-1}+\dots+(-1)^d dS_d&=dF \\
 &=-F\sum_{1\le i\le d}\frac{d{\it pr}^\ast_i t}{T-{\it pr}^\ast_i t} \\
 &=F\sum_{1\le i\le d}\frac{1}{{\it pr}^\ast_i t}\frac{d{\it pr}^\ast_i t}{1-\frac{T}{{\it pr}^\ast_i t}}\\
 &=F\sum_{r\ge 0}\omega_{r+1}T^r.
 \end{align*}
 Comparing the coefficients of $T^r$, we obtain equalities
 \begin{align*}
 S_d\omega_{1}&=\pm dS_d \\
 S_d\omega_2\pm S_{d-1}\omega_1 &=\pm dS_{d-1} \\
 & \vdots\\
 S_d\omega_{r+1}+(\text{a linear combination of }\omega_r,\dots,\omega_{r-d})&=0\ (r\ge d)\\
 & \vdots
 \end{align*}
 
 The assertion follows by induction on $r$.

For the case when $j<0$, take $F$ as $(1-{\it pr}^\ast_1 tT)\cdots(1-{\it pr}^\ast_d tT)$ and argue similarly.
 \qed}

 \medskip

\medskip
(Proof of Theorem \ref{Witt})

1. 
 Let $e_{R''/R'}$ be the ramification index of $R''/R'$. 
Let $e_{R''/R}$ be the ramification index of $R''/R$ induced by ${\it pr}_i$. By Lemma \ref{unif}, 
$e_{R''/R}$ is independent of $i$. From the definition of the filtrations, the map ${\it pr}^\ast_i:{\rm W}_{m+1}(K)\to{\rm W}_{m+1}(K'')$ 
sends ${\rm fil}_n{\rm W}_{m+1}(K)$ into ${\rm fil}_{ne_{R''/R}}{\rm W}_{m+1}(K'')$. 
Since $S_d$ is a uniformizer of $R'$ by Lemma \ref{unif}, the equality
\begin{equation*}
de_{R''/R}=e_{R''/R'}
\end{equation*}
holds. 
This shows the identity
\begin{equation*}
{\rm fil}_{\left\lfloor\frac{n}{d}\right\rfloor}{\rm W}_{m+1}(K')=
{\rm fil}_{ne_{R''/R}}{\rm W}_{m+1}(K'')\cap{\rm W}_{m+1}(K'),
\end{equation*}
hence the assertion.

2.
This follows from Lemma \ref{anbasis}.
\qed

 \section{Generalized jacobians and blow-ups of symmetric powers}\label{gen}
 
 Let $S$ be a scheme, $C$ be a projective smooth $S$-scheme whose geometric fibers are connected 
 and of dimension $1$. Let $\mathfrak{m}$ be an effective Cartier divisor of $C/S$, i.e. a closed subscheme of $C$ which is finite flat of finite presentation over $S$. We also call $\mathfrak{m}$ a modulus.  
Let ue denote, for $S$-schemes $T$, the projections $C\times_S T\to T$ by the same symbol ${\it pr}$.
 In this section, we recall and study the notion of generalized jacobian varieties. Let $d$ be an integer and $\mathfrak{m}$ be a modulus.  
 Let $T$ be an $S$-scheme. Consider a datum $(\mathcal{L},\psi)$ such that
 \begin{itemize}
 \item $\mathcal{L}$ is an invertible sheaf of $\deg=d$ on $C_T$.
 \item $\psi$ is an isomorphism $\mathcal{O}_{\mathfrak{m}_T}\to \mathcal{L}\lvert_{\mathfrak{m}_T}$.
 \end{itemize}
 We say that two such data $(\mathcal{L},\psi)$ and $(\mathcal{L'},\psi')$ are isomorphic 
 if there exists an isomorphism of invertible sheaves $f: \mathcal{L}\to \mathcal{L'}$ making 
 the following diagram commutes
 \begin{equation*}
 \xymatrix{
 \mathcal{O}_{\mathfrak{m}_T}\ar[rr]^{\psi}\ar[rd]_{\psi'}&&\mathcal{L}\lvert_{\mathfrak{m}_T}
 \ar[ld]^{f\rvert_{\mathfrak{m}_T}}\\
 &\mathcal{L'}\lvert_{\mathfrak{m}_T}.
 }
 \end{equation*}
For an $S$-scheme $T$, define a set 
   
 \begin{equation*}
 \Pic^{d,{\rm pre}}(T):=\left\{\text{the isomorphism class of }(\mathcal{L},\psi) \text{ defined as above}\right\}.
\end{equation*}
$\Pic^{d,{\rm pre}}$ extends in an obvious way to a presheaf on $\Sch/S$, which we denote by $\Pic^{d,{\rm pre}}$ also. Define $\Pic^d$ as the \'etale sheafification of $\Pic^{d,{\rm pre}}$.
Their fundamental properties which we use without proofs are:
\begin{itemize}
\item $\Pic^d$ are represented by $S$-schemes. When $\mathfrak{m}$ is faithfully flat over $S$, $\Pic^{d,{\rm pre}}$ are 
already \'etale sheaves. 

\item $\Pic^{0}$ is a smooth commutative group $S$-scheme with geometrically connected fibers.
\item $\Pic^d$ are $\Pic^0$-torsors.
\end{itemize}
$\Pic^0$ is called the generalized jacobian variety of $C$ with modulus $\mathfrak{m}$. When 
$\mathfrak{m}=0$, this is the jacobian variety of $C$. In this case, we also denote $\Jac^d$ 
for $\Pic^d$. Let $\mathfrak{m}_1$ and $\mathfrak{m}_2$ be moduli such that $\mathfrak{m}_1\subset\mathfrak{m}_2$. There exists a natural map from $Pic_{C,\mathfrak{m}_2}^d$ to $Pic_{C,\mathfrak{m}_1}^d$, restricting $\psi$. Since $\mathfrak{m}_2$ is a finite 
$S$-scheme, this map is a surjection as a morphism of \'etale sheaves.

For a finite flat $S$-scheme of finite presentation $D$, define a presheaf $\mathcal{O}_{D}^\times$ on $\Sch/S$ by sending an $S$-scheme $T$ to the multiplicative group 
$\Gamma(T,\mathcal{O}_{D\times_S T}^\times)$, which is called the Weil restriction 
of $\mathbb{G}_{{\rm m},D}$ to $S$. This is an \'etale sheaf, and represented by a smooth group $S$-scheme. When $D=S$, this is $\mathbb{G}_{{\rm m},S}$. Define a map $\mathbb{G}_{{\rm m},S}\to \mathcal{O}_{D}^\times$ from the map of $S$-schemes $D\to S$. When $\deg D$ is strictly positive everywhere on $S$, this is an injection of \'etale sheaves.

Consider a map $\mathcal{O}_{\mathfrak{m}}^\times \to \Pic^0$ sending $s\in\mathcal{O}_{\mathfrak{m}}^\times$ to the pair $(\mathcal{O}_C,\mathcal{O}_{\mathfrak{m}}
\overset{s}\to\mathcal{O}_{\mathfrak{m}})$. The image of this map coincides 
with the kernel of the map $\Pic^0\to\Jac^0$, and the kernel of the map $\mathcal{O}_{\mathfrak{m}}^\times \to \Pic^0$ is the image of $\mathbb{G}_{{\rm m},S}\to \mathcal{O}_{\mathfrak{m}}^\times$ induced by 
the morphism of $S$-schemes $\mathfrak{m}\to S$. In summary, if $\deg \mathfrak{m}$ is everywhere strictly positive, we have a short exact sequence:
\begin{equation*}
0\to\mathcal{O}_{\mathfrak{m}}^\times/\mathbb{G}_{{\rm m},S}\to\Pic^0\to\Jac^0\to 0.
\end{equation*}
In particular, when $C\to S$ has a section $P:S\to C$, ${\rm Pic}_{C,P}^0$ is isomorphic to $\Jac^0$. In this case, $\Jac^d$ has an expression as a sheaf which does not depend on the choice of $P$. Let $T$ be an $S$-scheme, and $\mathcal{L}_1$ and $\mathcal{L}_2$ are invertible sheaves of $\deg=d$ on $C_T$. Define an equivalence relation 
on $\Jac^{d,{\rm pre}}$ such that $\mathcal{L}_1$ and $\mathcal{L}_2$ are equivalent if and only if 
there exists an invertible sheaf $\mathcal{M}$ on $T$ such that $\mathcal{L}_1\cong\mathcal{L}_2
\otimes{\it pr}^\ast\mathcal{M}$. If $C\to S$ has a section, the quotient presheaf of $\Jac^{d,{\rm pre}}$ by this equivalence relation is an \'etale sheaf and coincides with the \'etale sheafification of $\Jac^{d,{\rm pre}}$ via the natural surjection. In particular, the identity map $\Jac^d\to\Jac^d$ corresponds to an equivalence class of invertible sheaves on $C\times_S\Jac^d$. 
In this paper, we call this class the universal class of invertible sheaves of $\deg=d$. 

From now on we fix a modulus $\Tilde{\mathfrak{m}}$. We call a modulus $\mathfrak{m}$ a submodulus if $\mathfrak{m}\subset\Tilde{\mathfrak{m}}$ holds. 
Until the last paragraph, we treat 
the case when submoduli considered are everywhere strictly positive on $S$. Let $\mathfrak{m}$ be a submodulus which is everywhere strictly positive. Then, $\Pic^{d}$ has an explicit expression as a sheaf, as explained before. 

Denote the genus of $C$ by $g$. This is a locally constant function on $S$. We consider a 
condition on an integer $d$ as below:
\begin{equation}\label{d}
d\geq \max\{2g-1+\deg\Tilde{\mathfrak{m}},\deg\Tilde{\mathfrak{m}}\}.
\end{equation}
When $S$ is quasi-compact, such a $d$ always exists. For an integer $d$ and a submodulus $\mathfrak{m}$, denote $d_\mathfrak{m}:=d-\deg\Tilde{\mathfrak{m}}+\deg\mathfrak{m}$. If $d$ satisfies the condition (\ref{d}), $d_\mathfrak{m}$ satisfies the condition (\ref{d}) with $\Tilde{\mathfrak{m}}$ replaced by $\mathfrak{m}$. 

Fix an integer $d$ satisfying the condition (\ref{d}). Let $T$ be an $S$-scheme and $\mathcal{L}$ be an invertible sheaf of $\deg=d$ on $C_T$. For every usual point $t\in T$, 
$R^1{\it pr}_{\ast}(\mathcal{L}(-\Tilde{\mathfrak{m}})\rvert_{C_t})$
 and $R^1{\it pr}_{\ast}(\mathcal{L}\rvert_{C_t})$ are 
zero by Serre duality and a degree argument. In this case, ${\it pr}_{\ast}\mathcal{L}(-\Tilde{\mathfrak{m}})$ and ${\it pr}_{\ast}\mathcal{L}$ are locally free sheaves and their formations commute with 
any base change, i.e. for any morphism of $S$-schemes $f: T'\to T$, the base change 
morphisms $f^{\ast}{\it pr}_{\ast}\mathcal{L}\to {\it pr}_{\ast}f^{\ast}\mathcal{L}$ and 
 $f^{\ast}{\it pr}_{\ast}(\mathcal{L}(-\Tilde{\mathfrak{m}}))\to
 {\it pr}_{\ast}f^{\ast}(\mathcal{L}(-\Tilde{\mathfrak{m}}))$ are isomorphisms. Also $ R^1{\it pr}_{\ast}f^{\ast}\mathcal{L}$ 
 and $R^1 {\it pr}_{\ast}f^{\ast}(\mathcal{L}(-\Tilde{\mathfrak{m}}))$ are zero.

Let $\mathfrak{m}$ be a submodulus. 
In this section, we construct a following commutative diagram of smooth $S$-algebraic spaces:
\begin{equation*}
\xymatrix{
&\Tilde{C}^{(d_\mathfrak{m})}_\mathfrak{m}\ar[r]\ar@{^{(}->}[d]&\Pic^{d_\mathfrak{m}}\ar
@{}[ld]|{\Box}\ar@{^{(}->}[d]^{(\ref{4})}\\
X_\mathfrak{m}\ar[r]^{\cong}\ar[rd]&\mathbb{P}(\mathcal{E}_\mathfrak{m})\ar[r]\ar[d]^{(\ref{2})}&
P^{d_\mathfrak{m}}_\mathfrak{m}\ar[d]^{(\ref{1})}\\
&C^{(d_\mathfrak{m})}&\Jac^{d_\mathfrak{m}}.
}
\end{equation*}

Let $\mathcal{L}$ be an invertible sheaf on $C_T$ for an $S$-acheme $T$. Denote $\mathcal{L}/(\mathcal{L}(-\mathfrak{m}))$ by $\mathcal{L}_{\mathfrak{m}}$.

For an $S$-scheme $T$, consider a pair $(\mathcal{L},\phi)$ such that $\mathcal{L}$ is an invertible sheaf of $\deg=d_\mathfrak{m}$ on $C_T$ and 
$\phi$ is an injection $\mathcal{O}_T\to{\it pr}_\ast\mathcal{L}_\mathfrak{m}$ such that 
the quotient ${\it pr}_\ast\mathcal{L}_\mathfrak{m}/\mathcal{O}_T$ is locally free. Call such pairs 
$(\mathcal{L},\phi)$ and $(\mathcal{L}',\phi')$ isomorphic if there exists an isomorphism $f:\mathcal{L}\overset\cong\to\mathcal{L}'$ such that the following diagram commutes
\begin{equation*}
\xymatrix{
&\mathcal{O}_T\ar[rd]^{\phi'}\ar[ld]_{\phi}\\
{\it pr}_\ast\mathcal{L}_\mathfrak{m}\ar[rr]_{{\it pr}_\ast f}&&{\it pr}_\ast\mathcal{L}'_\mathfrak{m}.
}
\end{equation*}
Define $P^{d_\mathfrak{m}}_\mathfrak{m}(T)$ as the set of isomorphism classes of such pairs. This is an \'etale sheaf on $\Sch/S$. 
Define a map 
\begin{equation}\label{1}
P^{d_\mathfrak{m}}_\mathfrak{m}\to\Jac^{d_\mathfrak{m}}
\end{equation}
by forgetting $\phi$. 
Let $X$ be a scheme, and $\mathcal{F}$ be a locally free sheaf of finite rank on $X$. 
We use a contra-Grothendieck notation for a projective space. Thus the $X$-scheme 
$\mathbb{P}(\mathcal{F})$ parametrizes invertible subsheaves of $\mathcal{F}$.

\begin{lm}\label{P}
The sheaf $P^{d_\mathfrak{m}}_\mathfrak{m}$ is represented by a proper smooth $S$-algebraic space. 
Assume that $C\to S$ has a section. Let $\mathcal{L}'$ be a representative invertible sheaf of 
the universal class. Then, as sheaves on $\Sch/\Jac^{d_\mathfrak{m}}$,  $P^{d_\mathfrak{m}}_\mathfrak{m}$ is isomorphic to the projectivization 
$\mathbb{P}({\it pr}_\ast\mathcal{L}'_\mathfrak{m})$ of ${\it pr}_\ast\mathcal{L}'_\mathfrak{m}$. 
\end{lm}
\proof{
Since $C\to S$ has a section \'etale locally on $S$, it is enough to consider the case when 
$C(S)$ is not empty. In this case, $\Jac^{d_\mathfrak{m}}$ has an explicit expression as a sheaf, as 
explained before. 

 Via the map (\ref{1}), 
we regard $P^{d_\mathfrak{m}}_\mathfrak{m}$ as a sheaf on $\Sch/\Jac^{d_\mathfrak{m}}$. 
Fix a representative invertible sheaf $\mathcal{L}'$ of the universal class. 
Let $\mathcal{N}$ be an element of $\mathbb{P}({\it pr}_\ast\mathcal{L}'_\mathfrak{m})(T)$, 
where $T$ is a $\Jac^{d_\mathfrak{m}}$-scheme. Let $\phi:\mathcal{O}_T\to
{\it pr}_\ast((\mathcal{L}'\otimes{\it pr}^\ast\mathcal{N}^{-1})_\mathfrak{m})$ 
be a morphism obtained by tensoring the inclusion $\mathcal{N}\hookrightarrow{\it pr}_\ast\mathcal{L}'_\mathfrak{m}$ with $\mathcal{N}^{-1}$. Then, the correspondence $\mathcal{N}\mapsto(\mathcal{L}'\otimes{\it pr}^\ast\mathcal{N}^{-1},\phi)$ defines a morphism of sheaves on $\Sch/\Jac^{d_\mathfrak{m}}$, $\mathbb{P}({\it pr}_\ast\mathcal{L}'_\mathfrak{m})
\to P^{d_\mathfrak{m}}_\mathfrak{m}$. This is an isomorphism. Indeed, we can construct its inverse as follows. 
Let $T$ be a $\Jac^{d_\mathfrak{m}}$-scheme and $(\mathcal{L},\phi)$ be an element of $P^{d_\mathfrak{m}}_\mathfrak{m}(T)$. 
Let $a:T\to\Jac^{d_\mathfrak{m}}$ be the structure map. Then, there exists an invertible sheaf $\mathcal{N}$ 
on $T$ such that $\mathcal{L}\otimes{\it pr}^\ast\mathcal{N}$ is isomorphic to $a^\ast\mathcal{L}'$. Such an $\mathcal{N}$ is unique since $C\to S$ has a section. Then, $\mathcal{N}\overset{\phi\otimes\mathcal{N}}\to{\it pr}_\ast((\mathcal{L}\otimes
{\it pr}^\ast\mathcal{N})_\mathfrak{m})\overset\cong\to{\it pr}_\ast a^\ast\mathcal{L}'_\mathfrak{m}$ 
is an element of $\mathbb{P}({\it pr}_\ast\mathcal{L}'_\mathfrak{m})(T)$. 
\qed}

Let $(\mathcal{L},\phi)$ be the universal element on $P^{d_\mathfrak{m}}_\mathfrak{m}$. Define  $\mathcal{E}_\mathfrak{m}$ as the $\mathcal{O}_{P^{d_\mathfrak{m}}_\mathfrak{m}}$-
module fitting in the following cartesian diagram
\begin{equation}\label{y}
\xymatrix{{\it pr}_\ast(\mathcal{L}(-\mathfrak{m}))\ar[r]\ar[d]&0\ar[d]\\
\mathcal{E}_\mathfrak{m}\ar[r]\ar[d]&\mathcal{O}_{P^{d_\mathfrak{m}}_\mathfrak{m}}\ar[d]_{\phi}\\
{\it pr}_{\ast}\mathcal{L}\ar[r]_{p\ \ \ }
&{\it pr}_{\ast}\mathcal{L}_\mathfrak{m},
}
\end{equation}
where the bottom horizontal arrow is the pushforward of the quotient map. Since all the right arrows are locally split injections and $p:{\it pr}_\ast\mathcal{L}\to{\it pr}_\ast\mathcal{L}_\mathfrak{m}$ is a surjection of locally free sheaves, $\mathcal{E}_\mathfrak{m}$ is locally free of finite rank and all the left arrows are locally 
split injections.

Let $\mathbb{P}(\mathcal{E}_\mathfrak{m})$ be the projectivization of $\mathcal{E}_\mathfrak{m}$. As a sheaf on $\Sch/S$, $\mathbb{P}(\mathcal{E}_\mathfrak{m})$ 
parametrizes triples $(\mathcal{L},\phi,\mathcal{M})$ such that $(\mathcal{L},\phi)$ is an element of 
$P^{d_\mathfrak{m}}_\mathfrak{m}$ and $\mathcal{M}$ is an invertible subsheaf of ${\it pr}_\ast\mathcal{L}\oplus\mathcal{O}_T$ such that the following diagram commutes 
\begin{equation}\label{z}
\xymatrix{
\mathcal{M}\ar[r]\ar[d]&\mathcal{O}_T\ar[d]_{\phi}\\
{\it pr}_{\ast}\mathcal{L}\ar[r]_{p\ \ }
&{\it pr}_{\ast}\mathcal{L}_\mathfrak{m},
}
\end{equation}
where the left vertical arrow (resp. top horizontal arrow) is the composition of the inclusion 
$\mathcal{M}\hookrightarrow{\it pr}_\ast\mathcal{L}\oplus\mathcal{O}_T$ 
and the first (resp. second) projection. 
This is a proper smooth $S$-algebraic space. 
\begin{lm}\label{div}
The map  ${\it pr}_\ast(\mathcal{L}(-\mathfrak{m}))\to\mathcal{E}_\mathfrak{m}$ in (\ref{y}) induces a 
closed immersion $\mathbb{P}({\it pr}_\ast\mathcal{L}(-\mathfrak{m}))\allowbreak\hookrightarrow\mathbb{P}(\mathcal{E}_\mathfrak{m})$. The closed subspace $\mathbb{P}({\it pr}_\ast\mathcal{L}(-\mathfrak{m}))$ is a hyperplane bundle of $\mathbb{P}(\mathcal{E}_\mathfrak{m})$. 
\end{lm}
\proof{
 The assertion follows from the exact sequense
 \begin{equation*}
 0\to{\it pr}_\ast(\mathcal{L}(-\mathfrak{m}))\to\mathcal{E}_\mathfrak{m}
 \to\mathcal{O}_{P^{d_\mathfrak{m}}_\mathfrak{m}}\to 0.
 \end{equation*}
\qed}

As a subsheaf of $\mathbb{P}(\mathcal{E}_\mathfrak{m})$, $\mathbb{P}({\it pr}_\ast\mathcal{L}
(-\mathfrak{m}))$ parametrizes triples $(\mathcal{L},\phi,\mathcal{M})$ such that the first 
projection $\mathcal{M}\to{\it pr}_\ast\mathcal{L}$ factors through ${\it pr}_\ast\mathcal{L}(-
\mathfrak{m})$.

Let $T$ be an $S$-scheme and $(\mathcal{L},\phi,\mathcal{M})$ be an element of $\mathbb{P}
(\mathcal{E}_\mathfrak{m})(T)$. Since the arrow $\mathcal{E}_\mathfrak{m}\to
{\it pr}_\ast\mathcal{L}$ in (\ref{y}) is locally a split injection, the first projection $\mathcal{M}\to{\it pr}_\ast\mathcal{L}$ is injective and the cokernel is locally free. Since these hold after any base 
change $t\to T$ from the spectrum of a field, the map ${\it pr}^\ast\mathcal{M}_t\to\mathcal{L}_t$ 
is injective for a usual point $t$ of $T$. Thus $\mathcal{O}_{C_T}\to \mathcal{L}\otimes{\it pr}^\ast\mathcal{M}^{-1}$ defines an effective Cartier divisor. Since $\deg(\mathcal{L}^{-1}\otimes{\it pr}^\ast\mathcal{M})$ equals to $-d_\mathfrak{m}$, ${\rm Spec}(\mathcal{O}_{C_T}/(\mathcal{L}^{-1}\otimes{\it pr}^\ast\mathcal{M}))$ is finite flat of finite presentation of $\deg=d_\mathfrak{m}$ over $T$ by the Riemann-Roch formula.

Let $C^{(d_{\mathfrak{m}})}$ be the $d_\mathfrak{m}$th symmetric product of $C$, which parametrizes effective Cartier divisor of $\deg=d_\mathfrak{m}$ on $C$. Define a map 
\begin{equation}\label{2}
\mathbb{P}(\mathcal{E}_\mathfrak{m})\to C^{(d_\mathfrak{m})}
\end{equation}
 sending $(\mathcal{L},\phi,\mathcal{M})$ to 
 ${\rm Spec}(\mathcal{O}_{C_T}/(
\mathcal{L}^{-1}\otimes{\it pr}^\ast\mathcal{M}))\subset C_T$. 

Let $Z_0$ be the closed subscheme of $C^{(d_\mathfrak{m})}$ defined by the map $C^{(d-\deg\Tilde{\mathfrak{m}})}\to C^{(d_\mathfrak{m})}$, adding ${\mathfrak{m}}$. Let $X_\mathfrak{m}$ be the blow-up of $C^{(d_\mathfrak{m})}$ along $Z_0$. We define a map
 \begin{equation}\label{3}
 h:X_\mathfrak{m}\to\mathbb{P}(\mathcal{E}_\mathfrak{m})
 \end{equation}
   as follows.

Let $D$ be the universal effective Cartier divisor on $C^{(d_\mathfrak{m})}$. 
Denote $\mathcal{O}_{C\times_k C^{(d_\mathfrak{m})}}(D)$ by $\mathcal{O}_C(D)$ and $\mathcal{O}_C(D)\otimes\mathcal{O}_{\mathfrak{m}\times_k C^{(d_\mathfrak{m})}}$ by $\mathcal{O}_C(D)_{\mathfrak{m}}$ for short. The composition of the natural maps $\mathcal{O}_{C}\to\mathcal{O}_C(D)\to\mathcal{O}_C(D)_{\mathfrak{m}}$ defines a map of 
locally free sheaves $\mathcal{O}_{C^{(d_\mathfrak{m})}}\to{\it pr}_\ast(\mathcal{O}_C(D)_{\mathfrak{m}})$ on $C^{(d_\mathfrak{m})}$.
After a base change $T\to C^{(d_\mathfrak{m})}$, this map becomes 
zero if and only if $T\to C^{(d_\mathfrak{m})}$ 
factors through $Z_0$. Thus the image of the dual $({\it pr}_\ast\mathcal{O}_C(D)_{\mathfrak{m}})^\vee\to \mathcal{O}_{C^{(d_\mathfrak{m})}}$ of this map is the ideal $\mathcal{I}$ defining $Z_0$. Let $\mathcal{L}:=\mathcal{O}_{C\times_S X_\mathfrak{m}}(D)
\otimes{\it pr}^\ast(\mathcal{I}\mathcal{O}_{X_\mathfrak{m}})$. Define $\phi:\mathcal{O}_{X_\mathfrak{m}}\to {\it pr}_\ast(\mathcal{O}_{C\times_S X_\mathfrak{m}}(D)\otimes{\it pr}^\ast(\mathcal{I}\mathcal{O}_{X_\mathfrak{m}}))_\mathfrak{m}$ to be the 
morphism obtained from the map $(\mathcal{I}\mathcal{O}_{X_\mathfrak{m}})^{-1}\to{\it pr}_\ast\mathcal{O}_
{C\times_S X_\mathfrak{m}}(D)_\mathfrak{m}$ by tensoring $\mathcal{I}\mathcal{O}_{X_\mathfrak{m}}$. 
Let $\mathcal{I}\mathcal{O}_{X_\mathfrak{m}}\to
{\it pr}_\ast\mathcal{L}$ be the map induced from the natural inclusion $\mathcal{O}_{X_\mathfrak{m}}
\to{\it pr}_\ast\mathcal{O}_{C\times_S X_\mathfrak{m}}(D)$ by tensoring $\mathcal{I}\mathcal{O}_{X_\mathfrak{m}}$. This map and the natural inclusion  $\mathcal{I}\mathcal{O}_{X_\mathfrak{m}}\to\mathcal{O}_{X_\mathfrak{m}}$ make the sheaf  $\mathcal{I}\mathcal{O}_{X_\mathfrak{m}}$ into a subsheaf of ${\it pr}_\ast\mathcal{L}\oplus\mathcal{O}_{X_\mathfrak{m}}$, which makes the diagram (\ref{z}) commutes. 
The triple $(\mathcal{L},\phi,\mathcal{I}\mathcal{O}_{X_\mathfrak{m}})$ defines a morphism 
$h:X_\mathfrak{m}\to\mathbb{P}(\mathcal{E}_\mathfrak{m})$. From the construction, it is obvious that $h$ is 
a morphism over $C^{(d_\mathfrak{m})}$. 

\begin{lm}
\label{comp}

1. As a subsheaf of $\mathbb{P}(\mathcal{E}_\mathfrak{m})$, $\mathbb{P}(\mathcal{E}_\mathfrak{m})\times_{C^{(d_\mathfrak{m})}}Z_0$ parametrizes 
triples $(\mathcal{L},\phi,\mathcal{M})$ such that the second projection $\mathcal{M}\to
\mathcal{O}$ are zero. 
As closed subspaces of $\mathbb{P}(\mathcal{E}_\mathfrak{m})$,  $\mathbb{P}(\mathcal{E}_\mathfrak{m})\times_{C^{(d_\mathfrak{m})}}Z_0$ and $\mathbb{P}({\it pr}_\ast\mathcal{L}(-\mathfrak{m}))$ are equal.  
In particular, $\mathbb{P}(\mathcal{E}_\mathfrak{m})\times_{C^{(d_\mathfrak{m})}}Z_0$ 
is a smooth divisor of $\mathbb{P}(\mathcal{E}_\mathfrak{m})$. 

2. Let $V$ be the complement of $Z_0$ in $C^{(d_\mathfrak{m})}$. As a subsheaf of $\mathbb{P}(\mathcal{E}_\mathfrak{m})$, $\mathbb{P}(\mathcal{E}_\mathfrak{m})\times_{C^{(d_\mathfrak{m})}}V$ parametrizes triples 
$(\mathcal{L},\phi,\mathcal{M})$ such that the second projection $\mathcal{M}\to \mathcal{O}$ 
is an isomorphism. The projection $\mathbb{P}(\mathcal{E}_\mathfrak{m})\times_{C^{(d_\mathfrak{m})}}V\to V$
is an isomorphism and its inverse  coincides with the restriction of $h$ to $V$.

\end{lm}
\proof{
We are considering the following diagram:
\begin{equation*}
\xymatrix{
\mathbb{P}(\mathcal{E}_\mathfrak{m})\times_{C^{(d_\mathfrak{m})}}Z_0\ar[r]\ar[d]&
\mathbb{P}(\mathcal{E}_\mathfrak{m})\ar[d]&\mathbb{P}(\mathcal{E}_\mathfrak{m})\times_{
C^{(d_\mathfrak{m})}}V\ar[l]\ar[d]\\
Z_0\ar[r]&C^{(d_\mathfrak{m})}&V.\ar[l]
}
\end{equation*}
1.
Let $(\mathcal{L},\phi,\mathcal{M})$ be an element of $\mathbb{P}(\mathcal{E}_\mathfrak{m})(T)$. This maps into $Z_0$ via the map $\mathbb{P}(\mathcal{E}_\mathfrak{m})\to C^{(d_\mathfrak{m})}$ if and only if the composition of  ${\it pr}^\ast\mathcal{M}\to\mathcal{L}\to\mathcal{L}_{\mathfrak{m}}$ is zero. Since the 
right vertical arrow of (\ref{z}) is an injection, this occurs if and only if the second projection $\mathcal{M}\to \mathcal{O}_T$ is zero. The second assertion is obvious from the definition and the expression of $\mathbb{P}
({\it pr}_\ast\mathcal{L}(-\mathfrak{m}))$ as a subsheaf. 
The last assertion is verified for  
$\mathbb{P}({\it pr}_\ast\mathcal{L}(-\mathfrak{m}))$ in Lemma \ref{div}.

2. 
Let $T$ be a $S$-scheme and $(\mathcal{L},\phi,\mathcal{M})$ be an element of $\mathbb{P}(\mathcal{E}_\mathfrak{m})(T)$. Let $t$ be a usual point of $T$. By 1., the pull back of the projection 
$\mathcal{M}\to\mathcal{O}_T$ by $t\hookrightarrow T$ is an isomorphism if and only if the image of $t$ by the map 
\begin{equation*}
T\overset{(\mathcal{L},\phi,\mathcal{M})}\to\mathbb{P}(\mathcal{E}_\mathfrak{m})\to C^{(d_\mathfrak{m})}
\end{equation*}
is in $V$.

Let $p:\mathbb{P}(\mathcal{E}_\mathfrak{m})\times_{C^{(d_\mathfrak{m})}}V\to V$ be the projection. 
Since $h:X_\mathfrak{m}\to\mathbb{P}(\mathcal{E}_\mathfrak{m})$ is a $C^{(d_\mathfrak{m})}$-morphism, $p\circ h\lvert_V$ is the identity. Let $(\mathcal{L},\phi,\mathcal{M})$ be an element of $\mathbb{P}(\mathcal{E}_\mathfrak{m})\times_{C^{(d_\mathfrak{m})}}V(T)$. 
Identify $\mathcal{M}$ and $\mathcal{O}_T$ by the second projection. By this rigidification, 
$(\mathcal{L},\phi,\mathcal{O}_T)$ is determined by the first projection. Thus $p$ is an injection as a morphism of sheaves. The assertion follows.
\qed}

After these preparations, we obtain the following:
\begin{thm}\label{isom}
The morphism $h:X_\mathfrak{m}\to\mathbb{P}(\mathcal{E}_\mathfrak{m})$ in  (\ref{3}) 
is an isomorphism.
\end{thm}
\proof{
By Lemma \ref{comp}.1, there exists a unique map $\mathbb{P}(\mathcal{E}_\mathfrak{m})\to X_\mathfrak{m}$ 
which is a lift of $\mathbb{P}(\mathcal{E}_\mathfrak{m})\to C^{(d_\mathfrak{m})}$. Since they are smooth separated over $S$ 
and contain a common $S$-dense open subset V  by Lemma \ref{comp}.2, the assertion follows. 
\qed}

Let $T$ be an $S$-scheme and $(\mathcal{L}, \psi)$ be an element of $\Pic^{d_\mathfrak{m}}(T)$. 
Define $\phi$ as the composition $\mathcal{O}_T\to{\it pr}_\ast\mathcal{O}_{\mathfrak{m}_T}
\overset{{\it pr}_\ast\psi}\to{\it pr}_\ast\mathcal{L}_\mathfrak{m}$. Then, the correspondence 
$(\mathcal{L}, \psi)\mapsto(\mathcal{L},\phi)$ defines a morphism 
\begin{equation}\label{4}
\Pic^{d_\mathfrak{m}}\to P^{d_\mathfrak{m}}_\mathfrak{m}.
\end{equation}

\begin{lm}\label{open}
The morphism $\Pic^{d_\mathfrak{m}}\to P^{d_\mathfrak{m}}_\mathfrak{m}$ in (\ref{4}) is an open immersion. 
The open subspace $\Pic^{d_\mathfrak{m}}$ parametrizes pairs $(\mathcal{L},\phi)$ such that 
the maps $\mathcal{O}_C\to\mathcal{L}_\mathfrak{m}$ obtained from $\phi$ by 
adjunction are surjective. 
\end{lm}
\proof{
It is obvious that this morphism is an injection of sheaves. Let $(\mathcal{L},\phi)$ be an 
element of $P^{d_\mathfrak{m}}_\mathfrak{m}(T)$. This element is in $\Pic^{d_\mathfrak{m}}$ if and only if the map $\mathcal{O}_{C_T}
\to\mathcal{L}_\mathfrak{m}$ obtained from $\phi$ by adjunction is a surjection. This is an 
open condition.
\qed}

\medskip

Next, we study behavior of various schemes (or algebraic spaces) 
when one replaces the modulus $\mathfrak{m}$. 
Let $\mathfrak{m}$ be a submodulus and  $\mathfrak{m}':=\Tilde{\mathfrak{m}}-\mathfrak{m}$.
Define a closed immersion $C^{(d-\deg\mathfrak{m}')}\to C^{(d)}$ by adding $\mathfrak{m}'$. We denote this closed subscheme of $C^{(d)}$ by $Z_\mathfrak{m}$. If $\mathfrak{m}_1
\subset\mathfrak{m}_2$, the inclusion $Z_{\mathfrak{m}_1}\subset
 Z_{\mathfrak{m}_2}$ holds. The closed immersion $Z_{\mathfrak{m}_1}\hookrightarrow Z_{\mathfrak{m}_2}$ is induced by adding $\mathfrak{m}_2-\mathfrak{m}_1$. This induces a map 
 \begin{equation}\label{5}
 X_{\mathfrak{m}_1}\hookrightarrow X_{\mathfrak{m}_2}
 \end{equation}
 of the blow-ups along $Z_0$. 
 Let $\mathfrak{m}_1$ and $\mathfrak{m}_2$ be submoduli such that $\mathfrak{m}_1\subset\mathfrak{m}_2$. 
Define a map $i_{\mathfrak{m}_1,\mathfrak{m}_2}:P^{d_{\mathfrak{m}_1}}_{\mathfrak{m}_1}\to P^{d_{\mathfrak{m}_2}}_{\mathfrak{m}_2}$ by sending $(\mathcal{L}_1,\phi_1)$ to 
$(\mathcal{L}_1(\mathfrak{m}_2-\mathfrak{m}_1),\phi)$, where $\phi$ is the composition of $\phi_1$ and 
the natural injection ${\it pr}_\ast(\mathcal{L}_1)_{\mathfrak{m}_1}\to{\it pr}_\ast\mathcal{L}_1(\mathfrak{m}_2-\mathfrak{m}_1)_{\mathfrak{m}_2}$. The map $i_{\mathfrak{m}_1,\mathfrak{m}_2}$ is a closed immersion. 
\begin{pr}\label{replace}

1. Let $\mathfrak{m}_1$ and $\mathfrak{m}_2$ be  submoduli such that $\mathfrak{m}_1\subset\mathfrak{m}_2$. 
As a subsheaf of $P^{d_{\mathfrak{m}_2}}
_{\mathfrak{m}_2}$, $P^{d_{\mathfrak{m}_1}}_{\mathfrak{m}_1}$ parametrizes pairs $(\mathcal{L}_2,\phi_2)$ such that the compositions 
$\mathcal{O}_C\overset{{\it pr}^\ast\phi_2}\to(\mathcal{L}_2)_{\mathfrak{m}_2}\to(\mathcal{L}_2)_{\mathfrak{m}_2
-\mathfrak{m}_1}$ are zero.  
The commutative diagram
\begin{equation*}
\xymatrix{
X_{\mathfrak{m}_1}\ar[r]\ar[d]&P^{d_{\mathfrak{m}_1}}_{\mathfrak{m}_1}\ar[d]\\
X_{\mathfrak{m}_2}\ar[r]&P^{d_{\mathfrak{m}_2}}_{\mathfrak{m}_2}
}
\end{equation*}
induced by (\ref{3}), (\ref{5}), and the projections $\mathbb{P}(\mathcal{E}_{\mathfrak{m}_i})\to 
P^{d_{\mathfrak{m}_i}}_{\mathfrak{m}_i}$  
is a cartesian diagram.

2. Assume that a submodulus $\mathfrak{m}$ is the sum $\sum_i \mathfrak{m}_i$ of 
submoduli of $\deg=1$. Let $\mathfrak{m}'_i:=\sum_{j\neq i}\mathfrak{m}_j$. Then, the open subspace $\Pic^{d_\mathfrak{m}}$ of $P^{d_\mathfrak{m}}_{\mathfrak{m}}$ is the complement of 
$P^{d_{\mathfrak{m}'_i}}_{\mathfrak{m}'_i}$ for all $i$. 
\end{pr}
\proof{
1. The first assertion is obvious from the definition of $i_{\mathfrak{m}_1,\mathfrak{m}_2}$. 
To prove the second assertion, it is enough to show that $\mathcal{E}_{\mathfrak{m}_1}\cong
i_{\mathfrak{m}_1,\mathfrak{m}_2}^\ast\mathcal{E}_{\mathfrak{m}_2}$ by Theorem \ref{isom}. 
Let $(\mathcal{L}_i,\phi_i)$ be the universal elements of $P^{d_{\mathfrak{m}_i}}_{\mathfrak{m}_i}$. The pull back of the cartesian diagram
\begin{equation*}
\xymatrix{
\mathcal{E}_{\mathfrak{m}_2}\ar[r]\ar[d]&\mathcal{O}_{P^{d_{\mathfrak{m}_2}}_{\mathfrak{m}_2}}\ar[d]
\\ {\it pr}_\ast\mathcal{L}_2\ar[r]& {\it pr}_\ast(\mathcal{L}_2)_{\mathfrak{m}_2}
}
\end{equation*}
by $i_{\mathfrak{m}_1,\mathfrak{m}_2}$ extends to the diagram
\begin{equation*}
\xymatrix{
i_{\mathfrak{m}_1,\mathfrak{m}_2}^\ast\mathcal{E}_{\mathfrak{m}_2}\ar[r]\ar[d]&
\mathcal{O}_{P^{d_{\mathfrak{m}_1}}_{\mathfrak{m}_1}}\ar[d]\\
 {\it pr}_\ast\mathcal{L}_1\ar[r]\ar@{^{(}->}[d]&
 {\it pr}_\ast(\mathcal{L}_1)_{\mathfrak{m}_1}\ar@{^{(}->}[d]\\
 {\it pr}_\ast(\mathcal{L}_1(\mathfrak{m}_2-\mathfrak{m}_1))\ar[r]&
 {\it pr}_\ast(\mathcal{L}_1(\mathfrak{m}_2-\mathfrak{m}_1))_{\mathfrak{m}_2},
}
\end{equation*}
where the two squares are cartesian diagrams, which shows the assertion. 

2. This follows from Lemma \ref{open} and 1.
\qed}

Define an $S$-scheme $\Tilde{C}^{(d_\mathfrak{m})}_\mathfrak{m}$ as the fibered product
\begin{equation}\label{cart}
\xymatrix{
\Tilde{C}^{(d_\mathfrak{m})}_\mathfrak{m}\ar[r]\ar[d]&\Pic^{d_\mathfrak{m}}\ar[d]\\
X_\mathfrak{m}\ar[r]&P^{d_\mathfrak{m}}_\mathfrak{m},
}
\end{equation}
 where the bottom horizontal map $X_\mathfrak{m}\to P^{d_\mathfrak{m}}_\mathfrak{m}$ is the 
composition $X_\mathfrak{m}\overset\cong\to\mathbb{P}(\mathcal{E}_\mathfrak{m})
\to P^{d_\mathfrak{m}}_\mathfrak{m}$. The $S$-scheme $\Tilde{C}^{(d_\mathfrak{m})}_\mathfrak{m}$ is a projective space bundle on $\Pic^{d_\mathfrak{m}}$. 
\begin{pr}\label{blowup}
The first projection $\Tilde{C}^{(d_\mathfrak{m})}_{\mathfrak{m}}\to X_\mathfrak{m}$ is an open immersion. Moreover, if $\mathfrak{m}$ is the sum 
$ \sum_{i}\mathfrak{m}_i$ of submoduli of $\deg=1$, 
$\Tilde{C}^{(d_\mathfrak{m})}_{\mathfrak{m}}$ coincides with the complement of $X_{\mathfrak{m}'_i}$ 
for all $i$, where $ \mathfrak{m}'_i:=\sum_{j\neq i}\mathfrak{m}_j$. 
\end{pr}
\proof{
These are consequences of Lemma \ref{open} and Corollary \ref{replace}.
\qed}

Let $U$ be the  complement of $\mathfrak{m}$ in $C$. 
Consider a map $U^{(d_\mathfrak{m})}\to\Pic^{d_\mathfrak{m}}$ sending $D$ to $(\mathcal{O}_C(D),\iota_D)$, where $\iota_d$ is the one induced from the natural identification $\mathcal{O}_{C\setminus D}\cong\mathcal{O}_C(D)\lvert_{C\setminus D}$. This makes the  following diagram commutes
\begin{equation}\label{symu}
\xymatrix{
U^{(d_\mathfrak{m})}\ar[r]\ar[d]&\Pic^{d_\mathfrak{m}}\ar[d]\\
X_\mathfrak{m}\ar[r]&P_\mathfrak{m}^{d_\mathfrak{m}},
}
\end{equation}
which induces an $X_\mathfrak{m}$-morphism $U^{(d_\mathfrak{m})}\to \Tilde{C}^{(d_\mathfrak{m})}_\mathfrak{m}$. This is an open immersion, since the vertical arrows of (\ref{cart}) and the left vertical arrow of (\ref{symu}) are open immersions. 
Combining the previous results, we obtain the following:
\begin{cor}\label{cor}
As an open subscheme of $\Tilde{C}^{(d_\mathfrak{m})}_\mathfrak{m}$, $U^{(d_\mathfrak{m})}$ is the complement of $\Tilde{C}^{(d_\mathfrak{m})}
\times_{C^{(d_\mathfrak{m})}} Z_0$. 
\end{cor}
\proof{
After a finite faithfully flat base change of $S$, we may assume that $\mathfrak{m}$ decomposes the sum $ \sum_{i}\mathfrak{m}_i$ of submoduli $\mathfrak{m}_i$ of $\deg=1$. 
Let $V$ be the complement of $Z_0$ in $C^{(d_\mathfrak{m})}$. Since $V\setminus U^{(d_\mathfrak{m})}$ is included 
in $ \cup_{i}Z_{\mathfrak{m}_i'}$, where  $ \mathfrak{m}'_i:=\sum_{j\neq i}\mathfrak{m}_j$, the assertion follows from Proposition \ref{blowup}.

\qed}

Now assume $\mathfrak{m}=0$. If $C$ has an $S$-valued point $P$, it is well-known that $C^{(d)}$ is a projective space bundle over 
$\Jac^{d}$ when $d\geq\max\{2g-1,0\}$, where $g$ is the genus of $C$. In other words, there exists a locally free sheaf 
$\mathcal{F}$ of finite rank on $\Jac^{d}$ such that $C^{(d)}$ is isomorphic to $\mathbb{P}(\mathcal{F})$. Classically this is proved using the Poincar\'e bundle. On the other hand, using Proposition \ref{blowup}, we might prove this fact with an extra condition $d\geq{\rm max}\{2g,1\}$, identifying ${\rm Pic}_{C,P}^d\cong\Jac^d$.

We have a following corollary:
\begin{cor}\label{pi1}
Assume that $S$ is connected noetherian. 
Let $\mathfrak{m}$ be a modulus $>0$ (resp. $=0$) of $C$ and $d$ be a sufficiently large  integer. Take a geometric point $\overline{x}$ on $\Tilde{C}^{(d)}_{\mathfrak{m}}$ (resp. on $C^{(d)}$) 
and denote $\overline{y}$ its image to $\Pic^d$. Then, the morphism of profinite groups 
$\pi_1(\Tilde{C}^{(d)}_{\mathfrak{m}},\overline{x})\to\pi_1(\Pic^d,\overline{y})$ (resp. $\pi_1(C^{(d)},\overline{x})\to\pi_1(\Pic^d,\overline{y})$) induced from the projection $\Tilde{C}^{(d)}_{\mathfrak{m}}\to
\Pic^d$ (resp. $C^{(d)}\to\Pic^d$) is an isomorphism. 
\end{cor}
\proof{
When $\mathfrak{m}=0$, let us also denote $\Tilde{C}^{(d)}_{\mathfrak{m}}$ for $C^{(d)}$. If $\mathfrak{m}>0$ 
(resp. $=0$), $\Tilde{C}^{(d)}_{\mathfrak{m}}$ is a projective space bundle over $\Pic^d$ (resp. after the base change from $S$ to an \'etale cover). In any case, the morphism $\Tilde{C}^{(d)}_{\mathfrak{m}}\to\Pic^d$ is proper surjective smooth with geometrically connected fibers. Take a geometric point $\overline{s}$ of $\Tilde{C}^{(d)}_{\mathfrak{m},\overline{y}}$ above $\overline{x}$.
Since the scheme $\Tilde{C}^{(d)}_{\mathfrak{m},\overline{y}}$ is simply connected, the homotopy exact sequence 
\begin{equation*}
\pi_1(\Tilde{C}^{(d)}_{\mathfrak{m},\overline{y}},\overline{s})\to\pi_1(\Tilde{C}^{(d)}_{\mathfrak{m}},\overline{x})\to
\pi_1(\Pic^d,\overline{y})\to1
\end{equation*}
implies the assertion.
\qed}

\section{Proofs}\label{pro}
In this section, we prove Theorem \ref{main} and Theorem \ref{main2}. 

First we prove some lemmas.

Let $C$ be a projective smooth geometrically connected curve over a perfect field $k$. Let $\mathfrak{m}$ be a modulus on $C$ and write $\mathfrak{m}=n_1 P_1+\dots+n_r P_r$, where $P_1,\dots,P_r$ are distinct closed points of $\mathfrak{m}$. Denote 
the complement of $\mathfrak{m}$ in $C$ by $U$. Let $d_i:=\deg P_i$. Take a positive integer $d$ 
so that $d\geq\deg\mathfrak{m}$.

\begin{lm}\label{l}
The morphism $\pi:C^{(n_1 d_1)}\times_k\dots\times_k C^{(n_r d_r)}
\times_k C^{(d-\deg\mathfrak{m})}\to C^{(d)}$, taking the sum, is \'etale at the 
generic point of the closed subvariety $\{n_1 P_1\}\times\dots\times\{n_r P_r\}\times C^{(d-\deg\mathfrak{m})}$ of $C^{(n_1 d_1)}\times_k\dots\times_k C^{(n_r d_r)}
\times_k C^{(d-\deg\mathfrak{m})}$.
\end{lm}
\proof{
We may assume that $k$ is algebraically closed (hence $d_i=1$ for all $i$). Since the map 
$\pi:C^{(n_1)}\times_k\dots\times_k C^{(n_r)}
\times_k C^{(d-\deg\mathfrak{m})}\to C^{(d)}$ is finite flat, it is enough to show that 
there exists a closed point $Q$ of $n_1P_1+\dots n_rP_r+C^{(d-\deg\mathfrak{m})}$ over which 
there are $\deg\pi$ points on $C^{(n_1)}\times_k\dots\times_k C^{(n_r)}
\times_k C^{(d-\deg\mathfrak{m})}$. Choose $Q$ as a point corresponding to a divisor 
$n_1P_1+\dots n_rP_r+P_{r+1}+\dots+P_{r+d-\deg\mathfrak{m}}$, where $P_1,\dots,P_{r+d-\deg\mathfrak{m}}$ are distinct points of $U(k)$. 

\qed}

\begin{lm}\label{pisur}
The morphism $\pi_1(U^{d})\to\pi_1(U^{(d)})$ induced from the natural 
projection $U^d\to U^{(d)}$ (base points are omitted) is surjective.
\end{lm}
\proof{
Since $U^d$ and $U^{(d)}$ are geometrically connected over $k$, 
it is enough to show the surjectivity after the base change to an algebraic closue $\bar{k}$ by considering the homotopy exact sequence $1\to\pi_1(U^d_{\bar{k}})\to\pi_1(U^d)\to\pi_1({\rm Spec}(k))\to 1$ and the counterpart of $U^{(d)}$. Assume that $k$ is algebraically closed. Take a closed point $P\in U^{(d)}$. Then, the fiber of $U^d\to U^{(d)}$ over the point $dP$ consists of one point 
$(P,P,\dots,P)$, which shows the surjectivity. 
\qed}
\medskip

 (Proof of Theorem \ref{main2})
 Let $C$ be a projective smooth geometrically connected curve over a perfect field $k$. 
 Let $\mathfrak{m}=n_1 P_1+\dots+n_r P_r\ (n_i\ge 1)$ be a modulus on $C$, and $U$ be its complement. 
   Set $A$ as the subgroup
of ${\rm H}^1(U,\mathbb{Q}/\mathbb{Z})$ 
consisting of a character $\chi$ such that $
{\rm Sw}_{P_i}(\chi)\leq n_i -1$ for $i=1,\dots,r$, 
and $B$ as the subgroup of  ${\rm H}^1(\Pic,\mathbb{Q}/\mathbb{Z})$ consisting of $\rho$ which is multiplicative. 

We construct a map  $\Psi:B\to A$. Take $\rho\in B$. Define $\chi$ to be the pull back of $\rho^1$ by the natural map $U\to\Pic^1$. We need to show that the ramification is bounded by $\mathfrak{m}$. Take a natural number $d$ large enough so that $d$ satisfies the condition (\ref{d}) for $\mathfrak{m}$. Consider the following 
commutative diagram
\begin{equation}\label{+}
\xymatrix{
U^d=
U\times\dots\times U\ar[rr]^{\pi}\ar[dd]&&U^{(d)}\ar[dd]^{p}\\\\
\Pic^1\times\dots\times\Pic^1\ar[rr]&&\Pic^d.
}
\end{equation}
By the multiplicativity of $\rho$, we know that $\pi^\ast p^\ast\rho^d=\chi^{\boxtimes d}$. Lemma 
\ref{pisur} implies that $p^\ast\rho^d=\chi^{(d)}$. We show that  ${\rm Sw}_{P_i}(\chi)\leq n_i-1$. 
We may asume that $k$ is algebraically closed (hence $d_i=1$). 
Corollary \ref{cor}, Lemma \ref{l}, and Lemma \ref{blprod} imply that the Swan conductor of 
$\chi^{(n_i)}$, with respect to the DVR at the generic point of the blow-up of $C^{(n_i)}$ along 
$n_iP_i$, is zero. 
By Theorem \ref{Witt}, we obtain ${\rm Sw}_{P_i}(\chi)\leq n_i-1$. Thus the map $B\to A$, pulling back by $U\to\Pic^1$, is well-defined. 
We denote this map by $\Psi$. 

First we show the injectivity of $\Psi$. Take $\rho$ from the kernel of $\Psi$. Since the multiplication map $\Pic^n\times\Pic^m\to\Pic^{n+m}$ and the two projections 
$\Pic^n\times\Pic^m\to\Pic^n,\Pic^m$ have geometrically connected fibers, the triviality of two of 
$\rho^n,\rho^m,\rho^{n+m}$ implies the triviality of the other. Thus it is enough to show the 
triviality of $\rho^d$ for sufficiently large $d$. 
Consider the diagram (\ref{+}). By Lemma \ref{pisur}, we know that $p^\ast\rho^d$ is trivial, which 
 implies that $\rho^d$ is trivial by Corollary \ref{pi1}. 
 
 The surjectivity of $\Psi$ is proved as follows. 
 Take $\chi\in A$. Let $d$ be an integer satisfying the condition (\ref{d}) for $\mathfrak{m}$. Proposition 
 \ref{blowup}, Theorem \ref{Witt}, and 
 Lemma \ref{l} imply that the character $\chi^{(d)}$ extends to a character $\Tilde{\chi}^{(d)}$ on 
 $\Tilde{C}^{(d)}_\mathfrak{m}$. Corollary \ref{pi1} implies that $\Tilde{\chi}^{(d)}$ descends to a character $\rho^d$ on $\Pic^d$. Let $d_1$ and $d_2$ be integers which satisfy the condition (\ref{d}). The commutative diagram 
 \begin{equation*}
 \xymatrix{
 U^{(d_1)}\times U^{(d_2)}\ar[d]\ar[r]&U^{(d_1+d_2)}\ar[d]\\
 \Pic^{d_1}\times\Pic^{d_2}\ar[r]^{\ \ \ q}&\Pic^{d_1+d_2}
 }
 \end{equation*}
 and the fact that the left vertical map has geometrically connected fibers show $q^\ast\rho^{d_1+d_2}=
 \rho^{d_1}\boxtimes 1+1\boxtimes\rho^{d_2}$. 
 Fix a non-zero effective Cartier divisor $D$ on $U$ such that $\deg D$ satisfies the condition (\ref{d}). Let $\xi$ be the pull back of $\rho^{\deg D}$ by the map 
 ${\rm Spec}(k)\to\Pic^{\deg D}$, corresponding to the point $D$. For an arbitrary integer $n$, take a natural number $m$ so large that the integer $n+m\deg D$ satisfies the condition (\ref{d}). Define $\rho^n:=f^\ast\rho^{n+m\deg D}\cdot a^\ast\xi^{-m\deg D}$, where $f:\Pic^n\to\Pic^{n+m\deg D}$ is the multiplication by $\mathcal{O}_C(mD)$ and $a:
 \Pic^n\to{\rm Spec}(k)$ is the structure map. This construction does not depend on $m$, since the 
 multiplicativity of $\rho^n$ is already verified for large $n$. By the same reason, the characters $\rho^n$ form a multiplicative character on $\Pic$. The equality $\chi=\Psi(\rho)$ follows from the commutative diagram
 \begin{equation*}
 \xymatrix{
 U\ar[r]^{({\rm id},g)\ \ \ \ \ \ }&U\times U^{(\deg D)}\ar[r]\ar[d]&U^{(\deg D+1)}\ar[d]\\
 &\Pic^1\times\Pic^{\deg D}\ar[r]&\Pic^{\deg D+1},
 }
 \end{equation*}
 where $g$ is the composition of the structure map $U\to{\rm Spec}(k)$ and the map 
 ${\rm Spec}(k)\to U^{(\deg D)}$ corresponding to the divisor $D$. 
 Indeed, the pullback of $\rho^{\deg D+1}$ by the map $U\to U\times U^{(\deg D)}\to U^{(\deg D+1)}
 \to \Pic^{\deg D+1}$ is $\chi\cdot b^\ast\xi$, where $b:U\to{\rm Spec}(k)$ is the structure map. 
 On the other hand, the pull back of $\rho^{\deg D+1}$ by the other way is $\Psi(\rho)\cdot b^\ast\xi$.

\qed

(Proof of Theorem \ref{main})
 Let $(G^0,G^1)$ be a connected abelian covering of $(\Pic^0,\Pic^1)$. Since the $d$th power of $\Pic^1$ is isomorphic to $\Pic^d$ as $\Pic^0$-torsors, the 
$d$th power $G^d$ of $G^1$ is naturally equipped with a compatible morphism $G^d\to\Pic^d$ of 
torsors.
Let $K$ be the kernel of the map $G^0\to\Pic^0$. This is a finite constant group since 
$G^0\to\Pic^0$ is a Galois isogeny. Take a non-trivial homomorphism  $\chi:K\to\mathbb{Q}/\mathbb{Z}$. This defines characters $\rho^d\in{\rm H}^1(\Pic^d,\mathbb{Q}/\mathbb{Z})$ for all $d$. From the construction, they form a multiplicative character on $\Pic$. Theorem \ref{main2} implies that the pull back 
of $\rho^1$ by $U\to\Pic^1$ is non-trivial and its ramification is bounded by $\mathfrak{m}$, which shows the first part of Theorem \ref{main}. 

Define the category $\mathcal{C}_1$ as the category of geometrically 
connected abelian coverings of $U$ whose ramifications are bounded by $\mathfrak{m}$ and the category $\mathcal{C}_2$ as the category of 
connected abelian coverings of $(\Pic^0,\Pic^1)$. We have 
constructed a functor $\Phi:\mathcal{C}_2\to\mathcal{C}_1$. We show that this functor is an equivalence of categories. We only treat the case when $k$ is algebraically closed. General case 
follows from this special case by using an argument of Galois descent. 
When $k$ is algebraically closed, 
the notion of $G^1$ is superfluous, i.e. if a closed point $P$ of $U$ is fixed, $\mathcal{C}_2$ is isomorphic to the category of abelian isogenies $G^0\to\Pic^0$ with $G^0$ is connected smooth algebraic group (cf. Theorem \ref{main3}). Let $\mathcal{C}_3$ be the last 
category. 

Fix a $k$-valued point $P$ of $U$. We show that the functor $\Phi':\mathcal{C}_3\to\mathcal{C}_1$,  pulling back by the 
morphism $U\to\Pic^0$, sending $Q\mapsto\mathcal{O}_C(Q-P)$, is an equivalence. 
The faithfulness is obvious since there only occur connected coverings. Let $G_1^0,G_2^0$ be 
elements of $\mathcal{C}_3$ and $V_1:=\Phi'(G_1^0),V_2:=\Phi'(G_2^0)$. The map ${\rm Hom}(G_i^0,G_i^0)\to{\rm Hom}(V_i,V_i)$ 
is bijective since the both have the same number of elements (and by the injectivity). 
Thus, to prove 
the fullness of $\Phi'$, it is enough to show that, if there is a map $V_1\to V_2$, there is a map 
$G_1^0\to G_2^0$. Let $K_i$ be the kernel of $G_i^0\to\Pic^0$. $K_i$ is canonically identified 
with the Galois group of $V_i\to U$. If there is a map $V_1\to V_2$, there is a map of abelian groups 
$h:K_1\to K_2$, which is independent of the choice of $V_1\to V_2$. 
We show that the commutativity of the following diagram,
\begin{equation}\label{K}
\xymatrix{
&\pi_1(\Pic^0)\ar[rd]_{p_2}\ar[ld]^{p_1}&\\
K_1\ar[rr]^h&&K_2
}
\end{equation}
where the maps from the top to the bottoms are the canonical surjections. Assume that there is an 
element $\sigma\in\pi_1(\Pic^0)$ such that $p_2(\sigma)\neq hp_1(\sigma)$. Take a group 
homomorphism $\rho^0:K_2\to\mathbb{Q}/\mathbb{Z}$ such that the images of $p_2(\sigma)$ and 
$hp_1(\sigma)$ are different. Since the characters $\rho^0 p_2$ and $\rho^0 hp_1$ are 
multiplicative and are pulled back to the same character via the map $U\to\Pic^0$, they are the same character, a contradiction. Thus the diagram (\ref{K}) is commutative, which implies that the quotient 
group $G_1^0/{\rm ker}(h)$ of $G_1^0$ is isomorphic to $G_2^0$. 

For the essential surjectivity, we argue as follows. Let $V$ be a connected cyclic covering of $U$. Take a character on $U$ whose kernel corresponds to $V$. By Theorem \ref{main3}, this character is a pull back of a multiplicative character $\rho^0$ on $\Pic^0$. Let $G^0$ be an \'etale covering of $\Pic^0$ corresponding to the kernel of $\rho^0$. We need to show that $G^0$ has a group structure. By the definition, $G^0$ is connected. From the multiplicativity of $\rho^0$, we know that there is a commutative 
diagram 
\begin{equation*}
\xymatrix{
G^0\times G^0\ar[r]^{m_G}\ar[d]&G^0\ar[d]\\
\Pic^0\times\Pic^0\ar[r]&\Pic^0.
}
\end{equation*}
Let us denote the map $m_G$ multiplicatively.
Let $F$ be the fiber of $G^0\to\Pic^0$ over $1\in\Pic^0$. For distinct points $y_1,y_2\in F$, 
the multiplication from right by $y_1$ and $y_2$, $G^0\to G^0$ are distinct. Indeed, Assume that $xy_1=xy_2$ for all $x\in G^0$. The multiplication from left by $x$, $G^0\to G^0$ is a $\Pic^0$-morphism and sends $y_1$ and $y_2$ to the same point, 
which implies that $y_1=y_2$ since $G^0$ is a connected covering of $\Pic^0$. 

Thus there exists an element $e\in F$ such that $xe=x$ for all $x\in G^0$. Next we show the commutativity of $m_G$. This 
follows from the fact that $G^0\times G^0$ is a connected covering of $\Pic^0\times\Pic^0$ and 
that the maps $G^0\times G^0\to G^0$, $(x,y)\mapsto xy$ and $(x,y)\mapsto yx$ send $(e,e)$ to the 
same point $e$. The associativity is proved in a similar way. Therefore it is verified that $G^0$ has a 
commutative group structure such that $G^0\to\Pic^0$ is a group homomorphism, hence an abelian 
isogeny. 
It is easy to show that $G^0$ is pulled back to $V$. 
 For a general $V$, use the fact that $V$ is a connected component of 
the finite projective limit of cyclic connected coverings which are quotient of $V$.
\qed

\section*{Acknowledgement}

The author would like to express his sincere gratitude to his advisor Professor Takeshi Saito for 
a lot of advice and his patient encouragement.

\end{document}